\documentclass{amsart}
\usepackage{amsmath,amsopn,amssymb,amsfonts,latexsym,amsthm}

\newcommand{\la}{{\lambda}}

\newcommand{\snknk}{{(S_k \times S_{n-k})\times (S_k \times
S_{n-k})}}

\newtheorem{thm}{Theorem}[section]

\newtheorem{lem}[thm]{Lemma}

\newtheorem{cor}[thm]{Corollary}
\newtheorem{fac}[thm]{Fact}
\newtheorem{obs}[thm]{Observation}

\newtheorem{exa}[thm]{Example}
\newtheorem{defn}[thm]{Definition}

\newcommand{\eqdef}{:=}

\newtheorem{prop}[thm]{Proposition}
\newtheorem{conj}[thm]{Conjecture}
\newtheorem{clm}[thm]{Claim}
\newcommand{\sumlim}{\sum\limits}
\newcommand{\een}{\end{enumerate}}
\newcommand{\blem}{\begin{lem}}
\newcommand{\elem}{\end{lem}}
\newcommand{\bcl}{\begin{clm}}
\newcommand{\ecl}{\end{clm}}
\newcommand{\ethm}{\end{thm}}
\newcommand{\bpr}{\begin{prop}}
\newcommand{\epr}{\end{prop}}
\newcommand{\bco}{\begin{cor}}
\newcommand{\eco}{\end{cor}}
\newcommand{\bcon}{\begin{conj}}
\newcommand{\econ}{\end{conj}}
\newcommand{\bde}{\begin{defn}}
\newcommand{\ede}{\end{defn}}
\newcommand{\bex}{\begin{exa}}
\newcommand{\eexa}{\end{exa}}
\newcommand{\bobs}{\begin{obs}}
\newcommand{\eobs}{\end{obs}}
\newcommand{\bexe}{\begin{exe}}
\newcommand{\eexe}{\end{exe}}

\newcommand{\grn}{G_{r,n}}

\begin{document}

\title{Permutation representations on invertible matrices}
\author{Yona Cherniavsky, Eli Bagno }
 \address{Department of Mathematics and Statistics
 Bar-Ilan University
 Ramat-Gan, Israel 52900}
 \email{cherniy@math.biu.ac.il , bagnoe@math.biu.ac.il}

\maketitle
\begin{abstract}
We discuss permutation representations which are obtained by the
natural action of $S_n \times S_n$ on some special sets of
invertible matrices, defined by simple combinatorial attributes.
We decompose these representations into irreducibles. The
multiplicities involved have a nice combinatorial interpretation.
We also generalize known results on asymptotic behavior of the
conjugacy representation of $S_n.$
\end{abstract}
\section{Introduction}

The (0,1)-matrices have a wide variety of applications in
combinatorics as well as in computer science. A lot of research
had been devoted to this area. By considering the set of $n \times
n$ (0,1)-matrices as a boolean monoid and relating them to posets,
one can get interesting representations of $S_n$. Note that $S_n
\times S_n$ acts on matrices by permuting rows and columns. Some
aspects of the corresponding equivalence relation are treated in
\cite{I} and \cite{Li}.  A simultaneous lexicographic ordering of
the rows and the columns using this action is shown in \cite{MM}.

The above action of $S_n \times S_n$ gives rise to a permutation
representation of $S_n \times S_n$ on (0,1)-matrices. If we
diagonally embed $S_n$ in $S_n \times S_n$ we get a generalization
of the conjugacy representation of $S_n$.

Adin and Frumkin \cite{AF} showed that the conjugacy character of
the symmetric group is close, in some sense, to the regular
character of $S_n$. More precisely, the quotient of the norms of
the regular character and the conjugacy character as well as the
cosine of the angle between them tend to 1 when $n$ tends to
infinity. This implies that these representations have essentially
the same decompositions.

Roichman \cite{R} further points out a wide family of irreducible
representations of $S_n$ whose multiplicity in the conjugacy
representation is asymptotically equal to their dimension, i.e.
their multiplicity in the regular representation.

In this paper we use the action of $S_n \times S_n$ on the (0,1)-
matrices to define two families of representations on a family of
orbits of this action. The first family forms an interpolation
between the regular representation of $S_n \times S_n$ and the
'diagonal sum' of the irreducible representations of $S_n$:
$\bigoplus_{\lambda \vdash n}{S^{\lambda} \otimes S^{\lambda}}$.
The other family is a generalization of the conjugacy
representation of $S_n$. In both cases we calculate characters and
present the decomposition of these representations into
irreducibles. The second family of representations can be seen as
an extension of the results of \cite{AF} and \cite{R}.

The rest of this paper is organized as follows: In Section
\ref{pre} we give some preliminaries on permutation groups and
their conjugacy classes and on representations of $S_n$. In
Sections \ref{alpha} and \ref{beta} we define the actions of $S_n
\times S_n$ and of $S_n$ on invertible matrices. In Sections
\ref{alpha 0-1} and \ref{beta 0-1} we specialize them into the
case of (0,1)-matrices. We also discuss the two permutation
representation they give rise to and calculate their characters.
In Section \ref{asy beta} we show some asymptotic results
concerning the generalization of the conjugacy representation of
$S_n$. In Section \ref{decompose alpha} we decompose the
representation of $S_n \times S_n$ mentioned above into
irreducible representations. In Section \ref{color} we consider
the actions described above on color permutation groups instead of
(0,1)-matrices.

\section{Preliminaries} \label{pre}
\subsection{Permutation Groups}

$S_n$ is the group of all bijections from the set $\{1...n\}$ to
itself. Every $\pi \in S_{n}$ may be written in {\em disjoint
cycle form} usually omitting the 1-cycles of $\pi$. For example,
$\pi =365492187$ may also be written as $\pi =(9,7,1,3,5)(2,6)$.
Given $\pi , \tau \in S_{n}$ let $\pi \tau \eqdef \pi \circ \tau $
(composition of functions) so that, for example,
$(1,2)(2,3)=(1,2,3)$. Note that two permutations are conjugate in
$S_n$ if and only if they have the same cycle structure. In this
paper we write $\pi\sim\sigma$ if the permutations $\pi$ and
$\sigma$ are conjugate in $S_n$. We denote by $\hat{S}_n$ the set
of conjugacy classes of $S_n$ and by $C_\pi\leq S_n$ the
centralizer subgroup of the element $\pi\in S_n$. Let
$C(\pi)\subseteq S_n$ denote the conjugacy class of the element
$\pi\in S_n$. By $supp(\pi)$ we mean the set of digits which are
not fixed by $\pi$. An element $\pi\in S_n$ with $|supp(\pi)|=t$
can be considered as an element of $S_t$ and then $C_{\pi}^t$
denotes the centralizer subgroup of the element $\pi$ in $S_t$
while $C^t(\pi)$ denotes the conjugacy class of the element $\pi$
in $S_t$.
$\pi_k\pi_{n-k}$ denotes an element of $S_k\times S_{n-k}$ where $\pi_k\in S_k$ and $\pi_{n-k}\in S_{n-k}$.\\
$C^{k\times (n-k)}(\pi_k\pi_{n-k})$ denotes the conjugacy class of the element $\pi_k\pi_{n-k}$ in $S_k\times S_{n-k}$.\\

There is an obvious embedding of $S_n$ in $GL_n(\mathbb{F})$ where
is $\mathbb{F}$ is any field. Just think about a permutation
$\pi\in S_n$ as an $n\times n$ matrix obtained from the identity
matrix by permutations of the rows. More explicitly: for every
permutation $\pi \in S_n$ we identify $\pi$ with the matrix:

\begin{center}
$$[\pi]_{i,j} =  \left\{ \begin{array}{cc}
{1}      &  { i=\pi(j)} \\
{0 }          & \text{otherwise}
\end{array}
\right. $$
\end{center}

\subsection{Color permutation groups}

For later use, we define here the color permutation groups. For
$r,n \in \mathbb{N}$, let $G_{r,n}$ denote the group of all $n$ by
$n$ monomial matrices whose non-zero entries are complex $r$-th
roots of unity. This group can also be described as the wreath
product $C_r \wr S_n$ which is the semi-direct product $C_r^n
\rtimes S_n$, where $C_r^n$ is taken as the subgroup of all
diagonal matrices in $\grn$. For $r=1$, $\grn$ is just $S_n$ while
for $r=2$, $\grn=B_n$, the Weyl group of type $B$.

\subsection{Representations}
\subsubsection{Permutation representations}

 In this work we deal mainly with permutation representations.
 Given an action of a group $G$ on a set $M$, the appropriate
 representation space is the space spanned by the elements of $M$
 on which $G$ acts by linear extension. We list two well known facts about permutation
 representations.

\begin{fac}\label{permrep}
The character of the permutation representation calculated at some
$g \in G$ equals to the number of fixed points under $g$.
\end{fac}

\begin{fac}
 The multiplicity of the trivial representation in a given
permutation representation is equal to the number of orbits under
the corresponding action.
\end{fac}

An important example we will use extensively in this work is the
conjugacy representation which is the permutation representation
obtained by the action of the group on itself by conjugation.

\subsubsection{Representations of $S_n$}
Let $n$ be a nonnegative integer. A {\it partition} of $n$ is an
infinite sequence of nonnegative integers with finitely many
nonzero terms $\lambda=(\lambda_1,\lambda_2,\ldots)$, where
$\lambda_1\ge\lambda_2\ge\ldots$ and  $\sumlim_{i=1}^{\infty}
\lambda_i =n$.

The sum $\sum \la_i=n$ is called the {\it size} of $\lambda$,
denoted $|\la|$; write also $\lambda\vdash n$. The number of parts
of $\la$, $\ell(\la)$, is the maximal $j$ for which $\la_j>0$. The
unique partition of $n=0$ is the {\it empty partition}
$\emptyset=(0,0,\dots)$, which has length $\ell(\emptyset):=0$.
For a partition $\lambda=(\lambda_1,\ldots,\lambda_k,\ldots)$
define the {\it conjugate partition}
$\lambda'=(\lambda'_1,\dots,\lambda'_i,\ldots)$ by letting
$\lambda'_i$ be the number of parts of $\lambda$ that are $\ge i$
$(\forall i\ge 1)$.

A partition $\lambda=(\lambda_1,...,\lambda_k)$ may be viewed as
the subset
$$\{(i,j) \mid 1 \leq i \leq k, 1 \leq j \leq \lambda_i\}
\subseteq \mathbb{Z}^2,$$ the corresponding {\it Young diagram}.
Using this interpretation we may speak of the intersection
$\lambda \cap \mu$, the set difference $\lambda \setminus \mu$ and
the symmetric set difference $\lambda \triangle \mu$ of any two
partitions. Note that $|\lambda \triangle
\mu|=\sumlim_{k=1}^{\infty}{|\lambda_k-\mu_k|}$.

 It is well known that the irreducible
representations of $S_n$ are indexed by partitions of $n$ (See for
example \cite{Sag}) and the representations of $S_n \times S_n$
are indexed by pairs of partitions $(\lambda,\mu)$ where
$\lambda,\mu \vdash n$. For every two representations of $S_n$,
$\lambda$ and $\rho$, we denote by $m(\lambda,\rho)$ the
multiplicity of $\lambda$ in $\rho$. If we denote by
$\langle\,,\rangle$ the standard scalar product of characters of a
finite group $G$ i.e.
$$
\langle\chi_1,\chi_2\rangle={1\over |G|}\sum_{\pi\in
G}\chi_1(\pi)\overline{\chi_2(\pi)}$$ then
$m(\lambda,\rho)=\langle\chi_{\lambda},\chi_{\rho}\rangle$.

Similarly, $m\big((\lambda,\mu),\varphi\big)$ denotes the
multiplicity of the representation of $S_n\times S_n$
corresponding to the pair of
 partitions $(\lambda,\mu)$ , $\lambda\vdash n$ , $\mu\vdash n$ in the
 decomposition of $\varphi$, where $\varphi$ is any representation of $S_n\times S_n$.

We cite here for later use the branching rule for the
representations of $S_n$. We start with a definition needed to
state the branching rule.

\bde Let $\lambda \vdash n$ be a Young diagram. Then a {\it
corner} of $\lambda$ is a cell $(i,j) \in \lambda$ such whose
removal leaves leaves the Young diagram of a partition. Any
partition obtained by such a removal is denoted by $\lambda^-$.
\ede

 \bpr \cite{Sag}\label{branching rule}
If $\lambda \vdash n$ then $$S^{\lambda} \downarrow
^{S_n}_{S_{n-1}} \cong \bigoplus_{\lambda^-}{S^{\lambda^-}}.$$
\epr

\section{The action of $S_n \times S_n$ on invertible
matrices}\label{alpha}

\bde Let  $G$ be a subgroup of $S_n \times S_n$ and let
$\mathbb{F}$ be any field. We define an action of $G$ on the group
$GL_n(\mathbb{F})$ by
$$
(\pi,\sigma)\bullet A = \pi A
\sigma^{-1}\,\,\,\text{where}\,\,(\pi,\sigma)\in G\,\,\text{and}\,\,A\in GL_n(\mathbb F)
\eqno(1)
$$
\ede
 It is easy to
see that this really defines a group action.

In this work we deal only with the cases: $G=S_n \times S_n$ and
$G=\left(S_k\times S_{n-k}\right) \times \left(S_k\times
S_{n-k}\right)$.

%

\bde Let $M$ be a finite subset of $GL_n(\mathbb{F})$, invariant
under the action of $S_n \times S_n$ defined above. We denote by
$\alpha_M$ the permutation representation of $G$ obtained from the
action $(1)$ . In the sequel we identify the action $(1)$ with the
permutation representation $\alpha_M$ associated with it. \ede

\subsection{A generalization of the conjugacy representation of
$S_n$}\label{beta}

In this section we present a conjugacy representation of $S_n$ on
a subset $M$ of $GL_n(\mathbb{F})$.

\bde Denote by $\beta$ the permutation representation of $S_n$
obtained by the following action on $M$.
$$
\pi\circ A=(\pi,\pi)\bullet A=\pi A\pi^{-1}\eqno(2)
$$
\ede

The connection between $\alpha_M$ and $\beta_M$ is given by the
following easily seen claim:
\bcl
Consider the diagonal embedding of $S_n$ into $S_n\times S_n$. Then
$$
\beta_M=\alpha_M\downarrow^{S_n\times S_n}_{S_n}.\eqno\qed
$$
\ecl
\begin{thm}\label{2chars}
For every  finite set $M \subseteq GL_n(\mathbb{F})$ invariant
under the action $(1)$ of $S_n \times S_n$ defined above:
\item
If $\pi$ and $\sigma$ are conjugate in $S_n$ then
$$ \chi_{\alpha_M}\left(
(\pi,\sigma)\right)=\chi_{\alpha_M}\left(
(\pi,\pi)\right)=\chi_{\beta_M}(\pi)=\#\{A\in M\,|\,\pi A=A\pi\}\,.
$$
\item
If $\pi$ is not conjugate to $\sigma$ in $S_n$ then
$$
\chi_{\alpha_M}\left((\pi,\sigma)\right)=0\,.
$$
\end{thm}
\begin{proof} See Theorem 4.5 in \cite{CS}.
\end{proof}
 \bex \label{alphan0}
  Take $M=S_n$ (embedded in $GL_n(\mathbb{F})$ as
permutation matrices). In this case $\beta_M$ is just the
conjugacy representation of $S_n$ and a direct calculation shows
that for every $\pi \in S_n$:
$$\chi_{\beta_M}(\pi)=|C_\pi|=\chi_{\alpha_M} (\pi,\pi).$$

For every irreducible representation of $S_n$ corresponding to a
partition $\lambda$ one has:

$$
m(\lambda,\beta_M)=\sum_C \chi_\lambda(C)\,\,\text{where }\,\,C
\text{ runs over all conjugacy classes of}\,\,S_n.
$$
Moreover,

$$
m\big( (\lambda,\lambda),\alpha_M\big)=1,
$$
and $$ m\big((\lambda,\mu),\alpha_M\big)=0\,\,\text{
when}\,\,\lambda\neq\mu.
$$
This means that $\alpha_M\cong\bigoplus_{\lambda\vdash
n}S^{\lambda}\otimes S^{\lambda}$ where $S^{\lambda}$ is the
irreducible $S_n$ - module corresponding to $\lambda$. This fact
will be implied by substitution $k=0$ in Proposition
\ref{charalpha}. \eexa

\section{The action of $S_n\times S_n$ on (0,1)-matrices} \label{alpha 0-1}
 In this section we specialize the action (1) of $S_n \times S_n$
 defined in Section \ref{alpha} to (0,1)-matrices.
 Consider the group $G=GL_n(\mathbb{Z}_2)$. For every $A \in G$ denote by $o(A)$ the  number of nonzero
entries in $A$. One can associate with $A$ a pair of partitions of
$o(A)$ with $n$ parts $(\eta(A) ,\theta(A))$ where  $\eta(A)$
describes the distribution of nonzero entries in the rows of $A$
and $\theta(A)$ describes the same distribution for columns. For
example, if:
$$
A=\begin{pmatrix}
1 &0 &0 &0\\
1 &1 &1 &0\\
0 &0 &1 &0\\
1 &1 &1 &1
\end{pmatrix}
$$
then $\eta(A)=(4,3,1,1)\vdash 9$ and $\theta(A)=(3,3,2,1)\vdash 9$. \\

If we fix a pair of partitions $(\eta ,\theta)$ then the set of
matrices corresponding to $(\eta ,\theta)$ is closed under the
action (1), but this action is not necessarily transitive on such
a set, i.e. it can be decomposed into a union of several orbits.
For example, if
$$
A=\begin{pmatrix}
1 &0 &0 &1\\
0 &1 &1 &0\\
0 &1 &0 &0\\
1 &0 &0 &0
\end{pmatrix} \,\,\text{and}\,\,\, B=\begin{pmatrix}
1 &1 &0 &0\\
1 &0 &1 &0\\
0 &1 &0 &0\\
0 &0 &0 &1
\end{pmatrix}
$$
then $\eta(A)=\theta(A)=\eta(B)=\theta(B)=(2,2,1,1)$ but it can be
easily shown that $A\neq \pi B\sigma$ for any $\pi,\sigma\in S_4$.

We present now a family of subsets of $GL_n(\mathbb Z_2)$ which
will be proven shortly to be orbits of our action:

\
\bde
$$
H_n^0=\{A\in G\,|\, \eta(A)=\theta(A)=(1,1,1,\ldots,1)=1^n\}
$$
$$
H_n^1=\{A\in G\,|\,
\eta(A)=(n,1,1,\ldots,1),\theta(A)=(2,2,\ldots,2,1)=2^{n-1}\,1\}
$$
$$
H_n^2=\{A\in G\,|\,
\eta(A)=(n,n-1,1,\ldots,1),\theta(A)=(3,3,\ldots,3,2,1)=3^{n-2}\,2\,1\}
$$
$$
\cdots\cdots\cdots\cdots\cdots\cdots\cdots\cdots\cdots\cdots\cdots\cdots\cdots\cdots\cdots
$$
\begin{multline}
H_n^k=\{ A\in G\,|\, \eta(A)=(n,n-1,\ldots,n-(k-1),1,\ldots,1),\qquad\qquad\\
\qquad\theta(A)=(k+1,k+1,\ldots,k+1,k,k-1,\ldots,2,1)=(k+1)^{n-k}\,k\,(k-1)\,\ldots\,2\,1\}
\nonumber
\end{multline}
$$
\cdots\cdots\cdots\cdots\cdots\cdots\cdots\cdots\cdots\cdots\cdots\cdots\cdots\cdots\cdots
$$
$$
H_n^n=\{A\in \,|\,
\eta(A)=\theta(A)=(n,n-1,n-2,\ldots,n-(k-1),\ldots,3,2,1)\}
$$
\ede
Note that in the above example $A\in H_4^2$.\\
A few remarks on the sets $H_n^k$ are in order: First, note that
$|H_n^k|=n!(n)_k$. Secondly, note that $H_n^0$ is $S_n$, embedded
as permutation matrices. Also note that the set $H_n^0\cup H_n^1$
is closed under matrix multiplication and matrix inversion and is
actually isomorphic to the group
$S_{n+1}$. Another simple observation is that $H_n^n=H_n^{n-1}$.\\

In order to prove that the sets $H_n^k$ are transitive under the
action we need the following definition:

\bde
  Denote by $U_{n,k}$ the following binary $n\times n$ matrix : the upper
left $k\times k$ block is upper triangular with the upper triangle
filled by ones, the upper right $k\times (n-k)$ block is filled by
ones, the lower left $(n-k)\times k$ block is the zero matrix and
the lower right $(n-k)\times (n-k)$ block is the identity matrix
$I_{n-k}$. \ede

\bpr \label{representative}
 Each set $H_n^k$ is transitive under the action (1). More
explicitly, $H_n^k=\{\pi U_{n,k}\sigma\,|\,\pi,\sigma\in S_n\}$
\epr
\begin{proof}

 We will prove that $U_{n,k}$ is a
representative of $H_n^k$.

Take an arbitrary matrix $A \in H_n^k$. By the definition of
$H_n^k$, $A$ has a unique row which is filled by ones. Choose a
permutation matrix $\pi_1$ such that this row is the first row of
$\pi_1 A$. By the definition of $H_n^k$, $\pi_1 A$ has a unique
column of type $(1,0,\ldots,0)^T$, thus one can choose a
permutation matrix $\sigma_1$ such that this is the first column
of $\pi_1 A\sigma_1$. Now consider the $(n-1)\times (n-1)$ sub
matrix $A_1$ obtained from $\pi_1 A\sigma_1$ by deleting the first
row and the first column. It is easy to see that $A_1\in
H_{n-1}^{k-1}$ and we can repeat the process described above. Note
that the first row and the first column of $A$ will not be changed
since their entries numbered $2..n$ are identical. (only $1$'s in
the first row, only $0$'s in the first column. Continuing this way
we get
$$
\pi_t\pi_{t-1}\cdots \pi_1
A\sigma_1\cdots\sigma_{t-1}\sigma_t=U_{n,k}
$$
which proves that $H_n^k\subseteq\{\pi
U_{n,k}\sigma\,|\,\pi,\sigma\in S_n\}$. The other inclusion, i.e.
$\{\pi U_{n,k}\sigma\,|\,\pi,\sigma\in S_n\}\subseteq H_n^k$,
follows directly from the definitions of $U_{n,k}$ and $H_n^k$.
\end{proof}

In particular
$$
H_n^0=\Bigg\{\pi\begin{pmatrix}
1 &0 &\cdots &\cdots &0\\
0 &1 &0      &\cdots &0\\
\vdots &\ddots&\ddots&\ddots&\vdots\\
0 &\cdots &0 &1 &0\\
0 &\cdots &\cdots &0 &1
\end{pmatrix}\sigma\,|\,\pi,\sigma \in S_n\Bigg\}=S_n
$$
$$
H_n^1=\Bigg\{\pi\begin{pmatrix}
1 &1 &1 &\cdots &1\\
0 &1 &0      &\cdots &0\\
\vdots &\ddots&\ddots&\ddots&\vdots\\
0 &\cdots &0 &1 &0\\
0 &\cdots &\cdots &0 &1
\end{pmatrix}\sigma\,|\,\pi,\sigma \in S_n\Bigg\}
$$
$$
H_n^3=\Bigg\{\pi\begin{pmatrix}
1 &1 &1 &1 &\cdots&\cdots &1\\
0 &1 &1 &1 &\cdots     &\cdots &1\\
0 &0 &1 &1 &\cdots    &\cdots &1\\
0 &0 &0 &1 &0 &\cdots &0\\
\vdots &\ddots&\ddots&\ddots&\ddots&&\ddots\vdots\\
0 &\cdots &\cdots &\cdots &0 &1 &0\\
0 &\cdots &\cdots &\cdots &\cdots&0 &1
\end{pmatrix}\sigma\,|\,\pi,\sigma \in S_n\Bigg\}
$$
\medskip
$$
\cdots\cdots\cdots\cdots\cdots\cdots\cdots\cdots\cdots\cdots\cdots\cdots\cdots\cdots\cdots\cdots\cdots\cdots
$$
\medskip
$$
H_n^n=\Bigg\{\pi\begin{pmatrix}
1 &1 &1 &1 &\cdots &1\\
0 &1 &1 &1      &\cdots &1\\
0 &0 &1 &1     &\cdots &1\\
\vdots &\ddots&\ddots&\ddots&\ddots&\vdots\\
0 &\cdots &\cdots &0 &1 &1\\
0 &\cdots &\cdots &\cdots &0 &1
\end{pmatrix}\sigma\,|\,\pi,\sigma \in S_n\Bigg\}
$$

For the case $k=n$ the permutation representation $\alpha_{H_n^k}$
can be easily described:

\bpr \label{alphann}

The representation $\alpha_{H^n_n}$ is isomorphic to the regular
representation of $S_n\times S_n$. \epr \noindent

\begin{proof}

As was shown in the proof of Proposition \ref{representative},
$H^n_n=\{\pi U_{n,n}\sigma\,|\,\pi,\sigma\in S_n\}$, where
$U_{n,n}$ is the upper triangular matrix with all the upper
triangle filled by ones. Define a bijection $\varphi: H_n^n \to
S_n \times S_n$ by $\varphi(\pi U_{n,n}\sigma)=(\pi,\sigma^{-1})$.
Since each row (column) of $U_{n,n}$ has a different number of
$1$-s (from 1 to $n$), we have:
$\pi_1U_{n,n}\sigma_1=\pi_2U_{n,n}\sigma_2$ if and only if
$\pi_1=\pi_2$ and $\sigma_1=\sigma_2$. This means that the mapping
$\varphi$ is well defined and bijective. Now:
$$(\omega,\tau) \bullet(\pi U_{n,n}\sigma)=\omega \pi U_{n,n}\sigma
\tau^{-1} \stackrel{\varphi}{\mapsto} (\omega \pi,\tau
\sigma^{-1})=(\omega,\tau)(\pi,\sigma^{-1}).$$

Thus $\varphi$ is an isomorphism of $S_n \times S_n$ - modules
between $H_n^n$ and the (left) regular representation of $S_n
\times S_n$.

\end{proof}

\subsection{A natural mapping from $H_n^k$ onto $S_n$.}\label{map}

In this section we present a surjection between the representation
of $S_n \times S_n$ on $H_n^k$ to the representation of $S_n
\times S_n$ on $S_n$. We will use this mapping later when we
decompose the permutation representation $\alpha$ into
irreducibles representations.

We deal first with $H_n^n$.

\bde Define the mapping $ t\,:\,H^n_n\longrightarrow S_n$ by

$$(\pi U_{n,n}\sigma) \rightarrowtail  \pi\sigma.$$
\ede This mapping is well defined since by the definition of
$U_{n,n}$
$$
\pi U_{n,n}\sigma=\rho U_{n,n}\tau\,\,\Longleftrightarrow
\pi=\rho\,\,\textrm{and}\,\,\sigma=\tau.
$$

It is easy to see that this mapping is a surjection and, moreover,
$|t^{-1}(\pi)|=n!$ for any $\pi\in S_n$. (Indeed, any fixed
$\pi\in S_n$ can be represented in exactly $n!$ ways in the form
$\pi=xy$ since choosing $x\in S_n$ in $n!$ ways we must take
$y=x^{-1}\pi$).

\bpr

The mapping $t$ preserves the action $\alpha$ of $S_n\times S_n$
on $H_n^n$, i.e.
$$
t(\pi A\sigma)=\pi t(A)\sigma\,\,\textrm{for any}\,\,A\in H^n_n.
$$
\epr \noindent
\begin{proof}
Any $A\in H^n_n$ can be represented $A=\rho U_{n,n}\tau$ for some
$\rho,\tau\in S_n$. Thus by the definition of $t$ we have
$$
t(\pi A\sigma)=t(\pi\rho U_{n,n}\tau\sigma)=\pi\rho\tau\sigma=\pi t(\rho U_{n,n}\tau)\sigma= \pi t(A)\sigma.
$$
\end{proof}

Now, in a similar way, for every $n,1 \leq k \leq n$ we define to $T_{n,k}:H_n^k\longrightarrow S_n$ as follows:

\bde Define the mapping $T_{n,k}\,:\,H_n^k\longrightarrow S_n$ by
$T_{n,k}(\pi U_{n,k}\sigma)=\pi\sigma$.

\ede

Notice that unlike the case of the sets $H_n^n$, here we have to
show that the mapping $T_{n,k}$ is well defined since a matrix $A
\in H_n^k$ can be written in the form $\pi U_{n,k} \sigma$ in more
than one way. Indeed, $T_{n,k}$ is well defined since if $\pi
U_{n,k} \sigma =U_{n,k}$ then as the first $k$ rows (columns) of
$U_{n,k}$ are different, we have $\pi , \sigma \in \{e\} \times
S_{n-k}$. Thus we can write $\pi=\pi_{n-k}$ and
$\sigma=\sigma_{n-k}$ and it is easy to see that one must have:
$\pi \sigma=e$ which says that $T_{n,k}(\pi U_{n,k}
\sigma)=e=T_{n,k}(U_{n,k}).$

Just as was the case in $H_n^n$ we have that:

\bpr \label{T_{n,k}}

The mapping $T_{n,k}$ preserves the action $\alpha$ of $S_n\times S_n$
on $H_n^k$, i.e.
$$
T_{n,k}(\pi A\sigma)=\pi T_{n,k}(A)\sigma\,\,\textrm{for any}\,\,A\in H_n^k.\eqno\qed
$$
\epr
It is also clear from the definition that $T_{n,k}$ is onto and it is easy to see that
$|T_{n,k}^{-1}(\pi)|=k!{n\choose k}=(n )_k$.

\section {The representation $\beta_M$ for $M=H_n^k$.}\label{beta 0-1}

In \cite{F} it was proven that the conjugacy representation of
$S_n$ contains every irreducible representation of $S_n$ as a
constituent. The representation $\beta$ defined in Section
\ref{beta} is a type of a conjugacy representation of $S_n$ on
$H_{n}^k$.

\bpr Denote the conjugacy representation of $S_n$ by $\psi$. Then
every irreducible representation of $S_n$ is a constituent in
$\beta_{H_n^k}$. In other words
$$
m\left(\lambda,\beta_{H_n^k}\right)>0\qquad\textrm{for any}\quad\lambda\vdash n.
$$ where
 $m\big(\lambda,\beta_{H_n^k}\big)$ denotes the multiplicity of
 the irreducible representation corresponding to $\lambda$ in
 $\beta_{H_n^k}$.
\epr
\begin{proof}  Denote the conjugacy
representation of $S_n$ by $\psi$.
 In [F] it is shown that $m\left(\lambda,\psi\right)>0$ for any $\lambda\vdash n$.
 By proposition \ref{T_{n,k}}, $T_{n,k}$ commutes with the conjugacy
 representation of $S_n$ and thus it gives rise to an epimorphism
 from $\beta_{H_n^k}$ onto the conjugacy representation of $S_n$. Hence, by Schur Lemma,
$$
m\left(\lambda,\beta_{H_n^k}\right)\geq m\left(\lambda,\psi\right)>0.
$$
\end{proof}
We turn now to the calculation of the character of $\beta_{H_n^k}$. By the definition, we have:

$$
\chi_{\beta_{H_n^k}}(\pi)\,(=\chi_{\alpha_{H_n^k}}(\pi,\pi))\,=\#\{A\in H_n^k\,|\, \pi A=A\pi\}
$$

but we can achieve much more than that:

\bpr \label{charbeta}
$$
\chi_{\beta_{H_n^k}}(\pi)=|C_{\pi}|(n-|supp(\pi)|)_k=(n-|supp(\pi)|)_k\chi_{Conj}(\pi)
$$
where $\chi_{Conj}$ is the conjugacy character of $S_n$. \epr

\begin{proof}
Denote $t=|supp(\pi)|$. We can assume that $\pi$ is of the
following form:

$$\pi=\begin{pmatrix}
\pi_t &0\\
0 &I_{n-t}
\end{pmatrix}\,,
$$
where $\pi_t\in S_t$ without fixed points.

 We have to calculate the number of matrices $A \in
H_n^k$ which commute with $\pi$. We recall the definition of
$H_n^k$:
\begin{multline}
H_n^k=\{ A\in G\,|\, \eta(A)=(n,n-1,\ldots,n-(k-1),1,\ldots,1),\qquad\qquad\\
\qquad\theta(A)=(k+1,k+1,\ldots,k+1,k,k-1,\ldots,2,1)=(k+1)^{n-k}\,k\,(k-1)\,\ldots\,2\,1\}
\nonumber
\end{multline}

For every $A \in H_n^k$ denote by $\delta(A)$ the vector
 of the row sums of $A$ and by $\varepsilon(A)$ the
vector of the column sums of $A$.

For example:

If $$ A=\begin{pmatrix}
1 &0 &0 &0\\
1 &1 &1 &0\\
0 &0 &1 &0\\
1 &1 &1 &1
\end{pmatrix} \in H_4^2
$$
then $\delta(A)=(1,3,1,4)\vDash 9$ and $\varepsilon(A)=(3,2,3,1)\vDash 9$.

Note that every $A \in H_n^k$ has $k$ rows with row sums
 ranging from $n$ to $n-k+1$, these will be called
'long rows'. The other $n-k$ rows are monomial. Likely, $A$ has
$k$ columns with column sums ranging from $1$ to $k$, these will
be called 'short columns'. The other $n-k$ columns have $k+1$
$1$-s each.

In $\pi A$ only the first $t$ rows of $A$ are permuted while in
$A\pi$ only the first $t$ columns of $A$ are permuted.

Note also that for every $\pi \in S_n$ we have:
$\delta(A\pi)=\delta(A)$ and $\varepsilon(\pi A)=\varepsilon(A)$.

Since all of the 'long' rows of $A$ have different row sums, if one
of the first $t$ rows of $A$ is 'long' then $\delta(\pi A) \neq
\delta(A)=\delta(A \pi)$ and thus $A\pi \neq \pi A$. Hence we can
assume that all the 'long' rows in $A$ are located after the first
$t$ rows of $A$. This implies that the first $t$
rows of $A$ are monomial. By similar arguments, the 'short' columns are
located after the first $t$ columns.

We claim now that the upper right $t \times (n-t)$ block of $A$ is the zero matrix.
(Indeed, if $A_{i,j}=1$ for some $1 \leq i \leq t$ and $t+1 \leq j \leq n$ then for each $1 \leq i' \leq t$ with
$i' \neq i$ we have: $A_{i',j} \neq 1$ since $A$ is invertible. Now, in $\pi A$ this $1$ moves to another
place while in $A \pi$ it is left in its original position).

Note that we have now that the upper left $t \times t$ block of $A$ is a
permutation matrix which commutes with $\pi_t$ in $S_t$.

 We turn now to the lower left $(n-t) \times t$ block of $A$. This block has $k$ 'long'
 rows and $n-t-k$ monomial rows. Note that by the consideration we described above,
 the monomial rows contribute only $0$ -s to this block. On the other hand, since the 'short' columns
 are located after the first $t$ columns, for every $1 \leq j \leq t$ the column numbered $j$ has exactly
 $k+1$ $1$-s. But for each $j$ one of this $1$-s must be located at the upper left $t \times t$ block since this block is a permutation matrix, hence
 the 'long' rows contribute only $1$-s to the lower left $(n-t) \times t$ block.

 Finally, as can be easily seen, the lower right $(n-t) \times (n-t)$ block is an arbitrary matrix of $H_{n-t}^{k}$.

Now calculate:
\begin{multline}
\chi_{\beta_{H_n^k}}(\pi)=\#\{A\in H_n^k\,|\, \pi A=A\pi\}=\\
=|C^t_{\pi}||H_{n-t}^k|=|C^t_{\pi}|(n-t)!(n-t)_k=\\
=|C_{\pi}|(n-t)_k=(n-|supp(\pi)|)_k\chi_{Conj}(\pi)\nonumber
\end{multline}
\end{proof}
Note that
$\chi_{\alpha_{H_n^k}}(\pi,\pi)=\chi_{\beta_{H_n^k}}(\pi)=0$ when
$t=|supp(\pi)|>n-k$. On the other hand, if $n-t \geq k$ then
$\chi_{\alpha_{H_n^k}}(\pi,\pi)=\chi_{\beta_{H_n^k}}(\pi)\neq 0$.

We turn now to the calculation of the multiplicity of every
irreducible representation of $S_n$ in $\beta_{H_n^k}$.
 \bpr
 Let $\lambda\vdash n$.
$$
m\big(\lambda,\beta_{H_n^k}\big)=\sum_{C\in
\hat{S}_n}\chi_{\lambda}(C)(n-|supp(C)|)_k
$$
where $\hat{S}_n$ denotes the set of conjugacy classes of
$S_n$.
\epr
\begin{proof} As above $t=|supp(\pi)|$. Then
\begin{multline}
m\big(\lambda,\beta_{H_n^k}\big)=\langle\chi_{\lambda},\chi_{\beta_{H_n^k}}\rangle=\frac{1}{|S_n|}\sum_{\pi\in S_n}\chi_{\lambda}(\pi)\chi_{\beta_{H_n^k}}(\pi)=\\
=\frac{1}{n!}\sum_{\pi\in S_n}\chi_{\lambda}(\pi)|C_{\pi}|(n-t)=
\sum_{\pi\in S_n}\chi_{\lambda}(\pi)\frac{n!}{|C(\pi)|}(n-t)_k=\\
=\sum_{\pi\in S_n}\chi_{\lambda}(\pi)\frac{1}{|C(\pi)|}(n-t)_k=\\
=\sum_{C\in
\hat{S}_n}|C|\chi_{\lambda}(C)\frac{1}{|C|}(n-|supp(C)|)_k=
\sum_{C\in \hat{S}_n}\chi_{\lambda}(C)(n-|supp(C)|)_k\nonumber
\end{multline}
\end{proof}

\section{Asymptotic behavior of the representation
$\beta_{H_n^k}$.}\label{asy beta}
 In this section we generalize the results of
Roichman \cite{R}, Adin, and  Frumkin \cite{AF} concerning the
asymptotic behavior of the conjugacy representation of  $S_n$.
These two results imply that the conjugacy representation and the
regular representation of $S_n$ have essentially the same
decomposition. In our case, as we prove in this section, the
representation $\beta_{H_n^k}$ is essentially $(n)_k$ times the
regular representation of $S_n$. We start by citing the result
from \cite{R}.

\noindent{\bf Theorem R1}  Let $m(\lambda)$ be the multiplicity of
the irreducible representation $S^\lambda$ in the conjugacy
representation of $S_n$, and let $f^\lambda$ be the multiplicity
of $S^\lambda$ in the regular representation of $S_n$. Then for
any $0 < \varepsilon < 1$ there exist $0<\delta(\varepsilon)$ and
$N(\varepsilon)$ such that, for any partition $\lambda$ of
$n>N(\varepsilon) $ with max$\{ {\lambda_1\over
n}\,,\,{\lambda_1'\over n}\}\leq\delta(\varepsilon)$,
$$
1-\varepsilon<\frac{m(\lambda)}{f^\lambda}<1+\varepsilon.
$$

The following generalization of this theorem is straightforward:

 \bpr \label{asymbeta}
 For any $0 < \varepsilon < 1$ there exist $0<\delta(\varepsilon)$ and
 $N(\varepsilon)$ such that, for any partition $\lambda$ of $n>N(\varepsilon)
 $ with max$\{ {\lambda_1\over n}\,,\,{\lambda_1'\over
 n}\}\leq\delta(\varepsilon)$,and for any $k \leq n$
$$
1-\varepsilon<\frac{m(\lambda,\beta_{H^k_n})}{(n)_k f^\lambda}<1+\varepsilon.
$$
\epr \begin{proof} In [R] it is shown that for any $0 <
\varepsilon <1 $ there exist $0<\delta(\varepsilon)$ and
$N(\varepsilon)$ such that, for any partition $\lambda$ of
$n>N(\varepsilon) $ with max$\{ {\lambda_1\over
n}\,,\,{\lambda_1'\over n}\}\leq\delta(\varepsilon)$,
$$
|m(\lambda)-f^{\lambda}|=|\sum_{C\in \hat{S}_n}\chi_{\lambda}(C)-f^{\lambda}|=|\sum_{C\neq id}\chi_{\lambda}(C)|\leq \varepsilon f^{\lambda},
$$
which immediately implies Theorem R1. In our case we have the trivial observation $(n-|supp(C)|)_k\leq (n)_k$ which together with the above gives us
\begin{multline}
|m(\lambda,\beta_{H^k_n})-(n)_k f^{\lambda}|=|\sum_{C\in \hat{S}_n}\chi_{\lambda}(C)(n-|supp(C)|)_k-(n)_k f^{\lambda}|=\\
=|\sum_{C\neq id}\chi_{\lambda}(C)(n-|supp(C)|)_k|\leq (n)_k|\sum_{C\neq id}\chi_{\lambda}(C)|\leq (n)_k\varepsilon f^{\lambda}\nonumber
\end{multline}
and our claim is proved.
\end{proof}

 The following asymptotic result from [AF] can also be generalized for the characters  $\chi_{\beta_{H^k_n}}$.\\
\noindent{\bf Theorem AF} Let $\chi_R^{(n)}$ and
$\chi_{Conj}^{(n)}$ be the regular and the conjugacy characters of
$S_n$ respectively. Then
$$
\lim_{n\to\infty}\frac{\|\chi_R^{(n)}\|}{\|\chi_{Conj}^{(n)}\|}=1\,,
$$
$$
\lim_{n\to\infty}\frac{\langle\chi_R^{(n)},\chi_{Conj}^{(n)}\rangle}{\|\chi_R^{(n)}\|\cdot\|\chi_{Conj}^{(n)}\|}=1\,
$$ where $\|*\|$ denotes the norm with respect to the standard
scalar product of characters.

 Our generalization looks as follows:
\bpr In the notations of Theorem AF
$$
\lim_{n\to\infty}\frac{\|(n)_k\chi_R^{(n)}\|}{\|\chi_{\beta_{H^k_n}}\|}=1\,,
$$
$$
\lim_{n\to\infty}\frac{\langle
(n)_k\chi_R^{(n)},\chi_{\beta_{H^k_n}}\rangle}{\|(n)_k\chi_R^{(n)}\|\cdot\|\chi_{\beta_{H^k_n}}\|}=
\lim_{n\to\infty}\frac{\langle\chi_R^{(n)},\chi_{\beta_{H^k_n}}\rangle}{\|\chi_R^{(n)}\|\cdot\|\chi_{\beta_{H^k_n}}\|}=
1\,,
$$

where $k$ is bounded or tends to infinity remaining less than $n$.

\epr

\begin{proof}
Denote for every $0 \leq k \leq n$
$$
F_k(S_n)=\sum_{C\in
\hat{S}_n}\frac{\big((n-|supp(C)|)_k\big)^2}{|C|}\,.
$$
Now
\begin{multline}
\|\chi_{\beta_{H^k_n}}\|^2=\langle\chi_{\beta_{H^k_n}},\chi_{\beta_{H^k_n}}\rangle=\\
={1\over n!}\sum_{\pi\in S_n}(\chi_{\beta_{H^k_n}}(\pi))^2
={1\over n!}\sum_{\pi\in S_n}|C_{\pi}|^2 \left((n-|supp(\pi)|)_k\right)^2=\\
={1\over n!}\sum_{\pi\in S_n}{(n!)^2\over |C(\pi)|^2}\left((n-|supp(\pi)|)_k\right)^2=\\
=n!\sum_{\pi\in S_n}{\big((n-|supp(\pi)|)_k\big)^2\over |C(\pi)|^2}=\\
=n!\sum_{C\in \hat{S}_n}|C|{\big((n-|supp(C)|)_k\big)^2\over |C|^2}=n!F_k(S_n)\nonumber
\end{multline}
It is well known that for every character $\chi$ of $S_n$:
$$
\langle\chi_R^{(n)},\chi\rangle=\chi(e)\quad\text{and}\quad
\|\chi_R^{(n)}\|^2=\langle\chi_R^{(n)},\chi_R^{(n)}\rangle=n! ,$$
hence
$$
\|\chi_{\beta_{H^k_n}}\|^2=F_k(S_n)\|\chi_R^{(n)}\|^2\,.
$$

Moreover, we have:

\begin{multline}
\frac{\langle\chi_R^{(n)},\chi_{\beta_{H^k_n}}\rangle}{\|\chi_R^{(n)}\|\cdot\|\chi_{\beta_{H^k_n}}\|}=
\frac{\chi_{\beta_{H^k_n}}(e)}{\sqrt{n!}\sqrt{n!
F_k(S_n)}}=\\
=\frac{(n)_k n!}{\sqrt{n!}\sqrt{n!
F_k(S_n)}}=\frac{(n)_k}{\sqrt{F_k(S_n)}}\,\,\,.\nonumber
\end{multline}

It is obvious from the definition of $F_k(S_n)$ that:
$$
\left((n)_k\right)^2\leq
F_k(S_n)\leq\left((n)_k\right)^2\sum_{C\in \hat{S}_n}{1\over
|C|}\,\,.
$$
In [AF] it is proved that
$$
\lim_{n\to\infty}\sum_{C\in \hat{S}_n}{1\over |C|}=1
$$
which establishes our claim.
\end{proof}

\section {The representations $\alpha_M$ for
$M=H_n^k$}\label{decompose alpha}

In this section we deal with the representations $\alpha_{H_n^k}$
defined in Section \ref{alpha}. We use the branching rule and the
Frobenious reciprocity to decompose these representations into
irreducible representations of $S_n \times S_n$. As we have
already seen in example \ref{alphan0}, $\alpha_{H_n^0} \cong
\bigoplus_{\lambda \vdash n}{S^{\lambda} \otimes S^{\lambda}}$
while $\alpha_{H_n^n}$ is the regular representation of $S_n
\times S_n \cong \bigoplus_{\lambda,\rho \vdash
n}{f^{\lambda}f^{\rho}S^{\lambda} \otimes S^{\rho}}$ and thus
$\alpha_{H_n^k}$ can be seen as a type of an interpolation between
these two representations.

First, concerning the character of $\alpha_{H_n^k}$, by combining
Proposition \ref{charbeta} and Claim \ref{2chars} together we get:

\begin{center}
$$\chi_{\alpha_{H_n^k}}(\pi,\sigma)=\left\{\begin{array}{cc}
|C_{\pi}|(n-|supp(\pi)|)_k\,\,,     &  {\pi} \text{ and } {\sigma} \text{ are conjugate in } {S_n}\\
{0 }\,\,,          & \text{otherwise}
\end{array}
\right.
$$
\end{center}
We turn now to the direct calculation of the multiplicity of an
irreducible representation of $S_n \times S_n$ in
$\alpha_{H_n^k}$.

\bpr For any $n$ and any $0\leq k\leq n$
$$
m\left( (\lambda,\mu),\alpha_{H^k_n}\right)=\frac{1}{n!}\sum_{\pi\in S_n}\chi_{\lambda}(\pi)\chi_{\mu}(\pi)(n-|supp(\pi)|)_k
$$
\epr \label{charalpha}

\noindent \begin{proof} Recall that $t=|supp(\pi)|$.
\begin{multline}
m\left( (\lambda,\mu),\alpha_{H^k_n}\right)=\langle\chi_{(\lambda,\mu)},\chi_{\alpha_{H_n^k}}\rangle_{S_n\times S_n}=\\
=\frac{1}{|S_n\times S_n|}\sum_{(\pi,\sigma)\in S_n\times S_n}\chi_{(\lambda,\mu)}\left((\pi,\sigma)\right)\chi_{\alpha_{H^k_n}}\left((\pi,\sigma)\right)=\\
= \frac{1}{(n!)^2}\sum_{(\pi,\sigma)\,:\,\pi\sim\sigma}\chi_{(\lambda,\mu)}\left((\pi,\sigma)\right)\chi_{\alpha_{H^k_n}}\left((\pi,\sigma)\right)=\\
=\frac{1}{(n!)^2}\sum_{\pi\in S_n}\sum_{\sigma:\pi\sim\sigma}\chi_{\lambda}(\pi)\chi_{\mu}(\pi)\chi_{\alpha_{H^k_n}}\left((\pi,\pi)\right)=\\
=\frac{1}{(n!)^2}\sum_{\pi\in S_n}|C(\pi)|\chi_{\lambda}(\pi)\chi_{\mu}(\pi)\chi_{\alpha_{H^k_n}}\left((\pi,\pi)\right)=\\
=\frac{1}{(n!)^2}\sum_{\pi\in S_n}|C(\pi)|\chi_{\lambda}(\pi)\chi_{\mu}(\pi)|C_{\pi}|(n-t)_k=\\
=\frac{1}{(n!)^2}\sum_{\pi\in S_n}|C(\pi)|\chi_{\lambda}(\pi)\chi_{\mu}(\pi)\frac{n!}{|C(\pi)|}(n-t)_k=\\
=\frac{1}{n!}\sum_{\pi\in
S_n}\chi_{\lambda}(\pi)\chi_{\mu}(\pi)(n-|supp(\pi)|)_k\nonumber
\end{multline}
\end{proof}
\noindent The boundary cases $k=0$ and $k=n$ are discussed in
Example \ref{alphan0} and Proposition \ref{alphann} respectively.

\noindent We list now several simple corollaries from Proposition
\ref{charalpha}.

 \bco
$$
m\left( (\lambda,\mu),\alpha_{H^k_n}\right)=m\left( (\lambda',\mu'),\alpha_{H^k_n}\right),
$$
where $\lambda'$ is the partition conjugate to $\lambda$, and
$$
m\left( (\lambda,\mu),\alpha_{H^k_n}\right)=m\left( (\mu,\lambda),\alpha_{H^k_n}\right).
$$
\eco
\begin{proof}

The first statement is due to the well known fact that
$\chi_{\lambda'}(\pi)=sgn(\pi)\chi_{\lambda}(\pi)$. The second one
is obvious.
\end{proof}

Note that $ m\left(\left( (n),(n)\right),\alpha_{H^k_n}\right)=1$
since the action $\alpha$ of $S_n \times S_n$ is transitive on
$H_n^k$. This implies the following corollary:

\bco

$$\sum_{\pi\in S_n}(n-|supp(\pi)|)_k=n!\qquad\text{for all }\quad
0\leq k\leq n.
$$ \qed
\eco



\bco
$$
m\left(\left( (n),1^n\right),\alpha_{H^k_n}\right)={1\over n!}\sum_{\pi\in S_n}sgn(\pi)(n-|supp(\pi)|)_k=0
$$
\eco

\begin{proof}

Clearly
\begin{multline}
m\left(\left( (n),1^n\right),\alpha_{H^k_n}\right)={1\over n!}\sum_{\pi\in S_n}sgn(\pi)(n-|supp(\pi)|)_k<\\
<{1\over n!}\sum_{\pi\in S_n}(n-|supp(\pi)|)_k=m\left(\left(
(n),(n)\right),\alpha_{H^k_n}\right)=1\nonumber \end{multline} But
$m\left(\left( (n),1^n\right),\alpha_{H^k_n}\right)$ is a non
negative integer, therefore it must be zero.
\end{proof}

\subsection{A combinatorial view of $\alpha_{H_n^k}$}

In this section we present another approach to the representation
$\alpha_{H_n^k}$. This approach will give us a combinatorial view
on the multiplicity formulas we calculated in the last section.

\bde Define the following subset of $H^k_n$:
$$
W^k_n=\left\{\pi_k\pi_{n-k} U_{n,k}\sigma_k\sigma_{n-k}\,|\,\pi_k,\sigma_k\in S_k\,\,\text{and}\,\,\pi_{n-k},\sigma_{n-k}\in S_{n-k}\right\}.
$$
\ede

The set $W^k_n$ is the orbit of the matrix $U_{n,k}$ under the
action $\alpha$ restricted to the subgroup $(S_k\times
S_{n-k})\times (S_k\times S_{n-k})$.

\bde Denote by $\omega_{n,k}$ the permutation representation of
the group $(S_k\times S_{n-k})\times (S_k\times S_{n-k})$ on
$W^k_n$ corresponding to the action $\alpha$. \ede

\bcl\label{omega}
$$
\omega_{n,k}\cong R_k\otimes\left(\bigoplus_{\rho\vdash
n-k}S^{\rho}\otimes S^{\rho}\right)
$$
where $R_k$ is the regular representation of $S_k\times S_k$. \ecl

\begin{proof}
We can view the action $\alpha$ of $\snknk$ on $W_n^k$ as composed
of two independent actions. One of them is the action of $S_k
\times S_k$ on $H_k^k$ (this corresponds to the the upper left
block of $U_{n,k}$) and the other is an action of $S_{n-k} \times
S_{n-k}$ on $S_{n-k}$ (this corresponds to the lower right block
of $U_{n,k}$). The permutation representation corresponding to the
first action is actually the regular representation of $S_k \times
S_k$ (see Proposition \ref{alphann}) while the second one is
$\bigoplus_{\rho\vdash n-k}S^{\rho}\otimes S^{\rho}$ (see Example
\ref{alphan0}).

\end{proof}

This implies the following:

\bcl \label{charomega}
$$
\chi_{\omega_{n,k}}(\pi_k\pi_{n-k},\sigma_k\sigma_{n-k})=
\left\{\begin{array}{ll}
0 & \textrm{when $\pi_k\neq e$ or $\sigma_k\neq e$}\\
0 & \textrm{when $\pi_{n-k}$ is not conjugate to $\sigma_{n-k}$ in $S_{n-k}$}\\
(k!)^2|C^{n-k}_{\pi_{n-k}}| & \textrm{when $\pi_k=\sigma_k=e$ and $\pi_{n-k}\sim\sigma_{n-k}$ in $S_{n-k}$}\end{array}\right.
$$
\qed \ecl \noindent

We can use $\omega_{n,k}$ to get information of $\alpha_{n,k}$.

\bpr \label{induce}
$$
\alpha_{H^k_n}=\omega_{n,k}\uparrow^{S_n\times S_n}_{(S_k\times S_{n-k})\times (S_k\times S_{n-k})}
$$
\epr
\noindent

\begin{proof}
Write $G=S_n \times S_n$ and $H=(S_k \times S_{n-k}) \times (S_k
\times S_{n-k})$ and identify $G/H$ with a prescribed set of left
transversals of $H$ in $G$. We clearly have:

$$H_n^k=\{g \bullet U_{n,k} \mid g \in S_n \times S_n\}=$$
$$\{(\sigma h) \bullet U_{n,k} \mid \sigma \in G/H, h \in H\}=$$
$$\{\sigma \bullet (h \bullet U_{n,k}) \mid \sigma \in G/H,h \in
H\}=$$
$$\{\sigma \bullet W_n^k \mid \sigma \in G/H\}=$$
$$\biguplus_{\sigma \in G/H}{\sigma \bullet W_n^k}$$
where $\uplus$ denotes disjoint union.

This implies that
$$
\alpha_{H^k_n}=\bigoplus_{\sigma \in G/H}{
 \sigma \bullet span_{\mathbb{C}}{W_n^k}}=\omega_{n,k}\uparrow^{S_n\times S_n}_{(S_k\times S_{n-k})\times
(S_k\times S_{n-k})}
$$ as claimed.
\end{proof}

We use now the Frobenius reciprocity to obtain the multiplicity of
any irreducible representation of $S_n \times S_n$ in
$\alpha_{H_n^k}$.

\bpr  Let $0 \leq k \leq n$ and let $\lambda, \mu$ be partitions
of $n$. Then
$$
m\left(\left(\lambda,\mu\right),\alpha_{H^k_n}\right)=\langle\chi_{\lambda}\downarrow^{S_n}_{S_{n-k}},\chi_{\mu}\downarrow^{S_n}_{S_{n-k}}\rangle
$$ or in other words:
$$
\alpha_{H^k_n}=\bigoplus_{\lambda,\mu\vdash
n}\langle\chi_{\lambda}\downarrow^{S_n}_{S_{n-k}},\chi_{\mu}\downarrow^{S_n}_{S_{n-k}}\rangle
S^{\lambda}\otimes S^{\mu}\eqno.
$$

\epr

\begin{proof}
Recall that $\chi_{(\lambda,\mu)}$ is the character of the
irreducible representation $S^{\lambda}\otimes S^{\mu}$ of
$S_n\times S_n$. Then, by Frobenius reciprocity, Claim
\ref{charomega} and Proposition \ref{induce} we have
\begin{multline}
m\left(\left(\lambda,\mu\right),\alpha_{H^k_n}\right)=\langle\chi_{\alpha_{H^k_n}},\chi_{(\lambda,\mu)}\rangle=\\
=\langle\chi_{\omega_{n,k}}\uparrow^{S_n\times S_n}_{(S_k\times S_{n-k})\times (S_k\times S_{n-k})},\chi_{(\lambda,\mu)}\rangle=\\
=\langle\chi_{\omega_{n,k}},\chi_{(\lambda,\mu)}\downarrow^{S_n\times S_n}_{(S_k\times S_{n-k})\times (S_k\times S_{n-k})}\rangle=\\
=\frac{1}{(k!)^2((n-k)!)^2}\sum_{(\pi_k\pi_{n-k},\sigma_k\sigma_{n-k})}\chi_{\omega_{n,k}}(\pi_k\pi_{n-k},\sigma_k\sigma_{n-k})\chi_{\lambda}(\pi_k\pi_{n-k})\chi_{\mu}(\sigma_k\sigma_{n-k})=\\
=\frac{1}{(k!)^2((n-k)!)^2}\sum_{\pi_{n-k}\sim\sigma_{n-k}}\chi_{\omega_{n,k}}(\pi_{n-k},\pi_{n-k})\chi_{\lambda}(\pi_{n-k})\chi_{\mu}(\pi_{n-k})=\\
=\frac{1}{(k!)^2((n-k)!)^2}\sum_{\pi_{n-k}\in S_{n-k}}|C^{n-k}(\pi_{n-k})|(k!)^2|C^{n-k}_{\pi_{n-k}}|\chi_{\lambda}(\pi_{n-k})\chi_{\mu}(\pi_{n-k})=\\
=\frac{1}{(n-k)!}\sum_{\pi_{n-k}\in S_{n-k}}\chi_{\lambda}(\pi_{n-k})\chi_{\mu}(\pi_{n-k})=\\
=\langle\chi_{\lambda}\downarrow^{S_n}_{S_{n-k}},\chi_{\mu}\downarrow^{S_n}_{S_{n-k}}\rangle\nonumber
\end{multline}
\end{proof}

The number
$\langle\chi_{\lambda}\downarrow^{S_n}_{S_{n-k}},\chi_{\mu}\downarrow^{S_n}_{S_{n-k}}\rangle$
has a very nice combinatorial interpretation. It follows from the
branching rule that this is just the number of ways to delete $k$
boundary cells from the diagrams corresponding to the partitions
$\lambda$ and $\mu$ to get the same Young diagram of $n-k$ cells.
By the branching rule (see Proposition \ref{branching rule}) we
have thus:

\bcl
$$
\langle\chi_{\lambda}\downarrow^{S_n}_{S_{n-k}},\chi_{\mu}\downarrow^{S_n}_{S_{n-k}}\rangle=0\,\,\text{ when }\,\,|\lambda\bigtriangleup\mu|>2k
$$
and it does not vanish otherwise.
 \ecl

 \bco
$$
m\left(\left(\lambda,\mu\right),\alpha_{H^k_n}\right)=0\,\,\text{ when }\,\,|\lambda\bigtriangleup\mu|>2k
$$
and
$$
m\left(\left(\lambda,\mu\right),\alpha_{H^k_n}\right)\neq
0\,\,\text{ when }\,\,|\lambda\bigtriangleup\mu|\leq 2k.\eqno \qed
$$
\eco

\subsection{Some asymptotic results}

In this section we use the fact that
$\beta_{H^k_n}=\alpha_{H^k_n}\downarrow_{S_n}^{S_n\times S_n}$ to
obtain some asymptotic results.

Embed $S_n$ in $S_n \times S_n$ as the diagonal subgroup. For
$\lambda, \mu, \nu \vdash n$, denote

$$
\gamma_{\lambda\mu\nu}=\frac{1}{n!}\sum_{\pi\in
S_n}\chi_{\lambda}(\pi)\chi_{\mu}(\pi)\chi_{\nu}(\pi).
$$ We have the following:

\bcl For $\lambda\vdash n$
$$
S^{\lambda}\uparrow_{S_n}^{S_n\times S_n}\cong \bigoplus_{\mu,\nu\vdash n}\gamma_{\lambda\mu\nu}S^{\mu}\otimes S^{\nu}
$$
\ecl \noindent
\begin{proof}

Let $\mu, \nu \vdash n$. By the Frobenius reciprocity,
\begin{multline}
m\left(\left(\mu,\nu\right),S^{\lambda}\uparrow_{S_n}^{S_n\times S_n}\right)=
\langle\chi_{(\mu,\nu)},\chi_{\lambda}\uparrow_{S_n}^{S_n\times S_n}\rangle=\\
=\langle \chi_{(\mu,\nu)}\downarrow_{S_n}^{S_n\times S_n},
\chi_{\lambda} \rangle =\frac{1}{n!}\sum_{\pi\in
S_n}\chi_{\lambda}(\pi)\chi_{\mu}(\pi)\chi_{\nu}(\pi)=\gamma_{\lambda\mu\nu}\nonumber
\end{multline}
\end{proof}

\noindent \textbf{Remark: } The numbers $\gamma_{\lambda\mu\nu}$
appear in [Md, page 115] in the context of the Schur functions
within the following formula:
$$
s_{\lambda}(xy)=\sum_{\mu,\nu}\gamma_{\lambda\mu\nu}s_{\mu}(x)s_{\nu}(y),
$$
where $x=(x_1,x_2,\ldots)$, $y=(y_1,y_2,\ldots)$ and $(xy)$ means the set of variables $x_i y_j$.

\noindent We have now the following asymptotic result:

 \bpr For any
$0 < \varepsilon <1 $ there exist $0<\delta(\varepsilon)$ and
$N(\varepsilon)$ such that, for any partition $\lambda$ of
$n>N(\varepsilon) $ with max$\{ {\lambda_1\over
n}\,,\,{\lambda_1'\over n}\}\leq\delta(\varepsilon)$,
$$
1-\varepsilon<\frac{\sum_{\mu,\nu\vdash n}\langle\chi_{\mu}\downarrow^{S_n}_{S_{n-k}},\chi_{\nu}\downarrow^{S_n}_{S_{n-k}}\rangle\gamma_{\lambda\mu\nu}}{(n)_k f^\lambda}<1+\varepsilon.
$$
\epr
\noindent

\begin{proof}

By Proposition \ref{asymbeta} we have:
$$
1-\varepsilon<\frac{m\left(\lambda,\beta_{H^k_n}\right)}{(n)_k f^\lambda}<1+\varepsilon.
$$
so we have to show that
$$
m\left(\lambda,\beta_{H^k_n}\right)=\sum_{\mu,\nu\vdash n}\langle\chi_{\mu}\downarrow^{S_n}_{S_{n-k}},\chi_{\nu}\downarrow^{S_n}_{S_{n-k}}\rangle\gamma_{\lambda\mu\nu}.
$$
Indeed by the previous claims and the  Frobenius reciprocity:
\begin{multline}
m\left(\lambda,\beta_{H^k_n}\right)=\langle\chi_{\beta_{H^k_n}},\chi_{\lambda}\rangle=\\
=\langle\chi_{\alpha_{H^k_n}}\downarrow_{S_n}^{S_n\times S_n},\chi_{\lambda}\rangle=\langle\chi_{\alpha_{H^k_n}},\chi_{\lambda}\uparrow_{S_n}^{S_n\times S_n}\rangle=\\
=\sum_{\mu,\nu\vdash
n}\langle\chi_{\mu}\downarrow^{S_n}_{S_{n-k}},\chi_{\nu}\downarrow^{S_n}_{S_{n-k}}\rangle\gamma_{\lambda\mu\nu}.\nonumber
\end{multline}
\end{proof}

Substituting in the above proposition $k=0$ and $k=n$ we get the
following:

\bco

$$ 1-\varepsilon<\frac{\sum_{\mu\vdash
n}\gamma_{\lambda\mu\mu}}{f^\lambda}<1+\varepsilon,
$$
$$ 1-\varepsilon<\frac{\sum_{\mu,\nu\vdash
n}\gamma_{\lambda\mu\nu}f^{\mu}f^{\nu}}{n!
f^\lambda}<1+\varepsilon.\eqno\square
$$

\eco
\textbf{Remark: }The first statement of this corollary follows from Theorem R1 and the equality
$$
\sum_{\mu\vdash
n}\gamma_{\lambda\mu\mu}=\sum_{C\in\hat{S}_n}\chi_{\lambda}(C)
$$
which follows from the character orthogonality relation.

\section{The actions $\alpha$ and $\beta$ on colored
permutations}\label{color}
 In this section we introduce actions of
$S_n$ and $S_n \times S_n$ on another family of sets, namely the
colored permutation groups. We start with the actions on $B_n= C_2
\wr S_n$.

\subsection{The action $\alpha$ of $S_n\times S_n$ on signed
permutations.}


Consider the action $\alpha$ of $S_n \times S_n$ on $B_n$. We
start by describing the orbits of this action.

\bde For every $0 \leq  k \leq n$ define

$$X^{k}_{n}=\{A\in B_n\,|\,A \,\,\text{has  exactly
}\,k\,\,\,\text{minuses}\}.
$$
\ede

For example
$$
\begin{pmatrix}
0 & 0 & 0 & 1\\
-1 & 0 & 0 & 0\\
0 & 1 & 0 & 0\\
0 & 0 & -1 & 0
\end{pmatrix}\in X_4^2.
$$

It is easy to see that the sets $X_n^k$ form a partition of $B_n$.
Also, note that $|X^k_n|=n!{n\choose k}$.

\bcl Each set $X^{k}_{n}$ is an orbit under the action $\alpha$ of
$S_n\times S_n$ on $B_n$, i.e.
$$
X^{k}_{n}=\{\pi \tilde{U}_{n,k}\sigma\,|\,\pi,\sigma\in S_n\},
$$
where
$$
\tilde{U}_{n,k}=\begin{pmatrix}
-I_k &0_{k\times (n-k)}\\
0_{(n-k)\times k} &I_{n-k}\end{pmatrix}
$$
and $I_t$ is the identity $t\times t$ matrix. \ecl

\begin{proof}
Let $A \in X_n^k$. Note that due to the semi direct decomposition
$B_n= C_2^n \rtimes S_n$ we have $A=ZP$ where $Z$ is a diagonal $n
\times n$ matrix with only $\pm 1$-s on the main diagonal and $P$
is a permutation matrix. Assuming that the rows numbered
$i_1,...,i_k$ in $A$ have minuses, denoting $\tau=(1,i_1) (2,i_2)
\cdots (k,i_k)$ we have $\tau A P^{-1} \tau^{-1}=\tilde{U}_{n,k}$.

\end{proof}

We decompose now the representations $\alpha_{X_n^k}$ into
irreducible representations just as we did in the previous
section.

 \bde Define the following
subset of $X^k_n$:
$$
\tilde{W}^k_n=\left\{\pi_k\pi_{n-k} \tilde{U}_{n,k}\sigma_k\sigma_{n-k}\,|\,\pi_k,\sigma_k\in S_k\,\,\text{and}\,\,\pi_{n-k},\sigma_{n-k}\in S_{n-k}\right\}.
$$
\ede The set $\tilde{W}^k_n$ is the orbit of the matrix
$\tilde{U}_{n,k}$ under the action $\alpha$ by the group
$(S_k\times S_{n-k})\times (S_k\times S_{n-k})$. \bde Denote
$\tilde{\omega}_{n,k}$ the permutation representation of the group
$(S_k\times S_{n-k})\times (S_k\times S_{n-k})$ which is obtained
from the action $\alpha$ of this group on the set $\tilde{W}^k_n$.
\ede

The proof of the following simple observation is similar to the
proof of Claim \ref{omega}.

\bcl

\item $$ \tilde{\omega}_{n,k}\cong \left(\bigoplus_{\rho\vdash
k}S^{\rho}\otimes
S^{\rho}\right)\otimes\left(\bigoplus_{\rho\vdash
n-k}S^{\rho}\otimes S^{\rho}\right)
$$

$$\chi_{\tilde{\omega}_{n,k}}(\pi_k\pi_{n-k},\sigma_k\sigma_{n-k})=
\left\{\begin{array}{ll}
|C^k_{\pi_k}||C^{n-k}_{\pi_{n-k}}| & \textrm{when $\pi_k\pi_{n-k}\sim\sigma_k\sigma_{n-k}$ in $S_k\times S_{n-k}$}\\
0 & \textrm{otherwise}\end{array}\right.
$$

\ecl \qed

\bpr
$$
\alpha_{X^k_n}=\tilde{\omega}_{n,k}\uparrow^{S_n\times S_n}_{(S_k\times S_{n-k})\times (S_k\times S_{n-k})}
$$
\epr
\begin{proof}Very similar to the proof of Prop. {\ref{induce}}.
\end{proof}

Recall from \cite{Sag} the definition of $c^{\lambda}_{\rho\nu}$
-- the Littlewood-Richardson coefficients defined by the following
formula:

$$
\left(S^{\rho}\otimes S^{\nu}\right)\uparrow^{S_n}_{S_k\times S_{n-k}}=
\bigoplus_{\lambda\vdash n}c^{\lambda}_{\rho\nu}S^{\lambda},
$$
where $\rho \vdash k$ and $\nu\vdash n-k$. Using the Frobenius
reciprocity formula we have for every $\lambda \vdash n$:

\bcl \label{LR}
$$
S^{\lambda}\downarrow^{S_n}_{S_k\times
S_{n-k}}=\bigoplus_{\rho\vdash k,\nu\vdash
n-k}c^{\lambda}_{\rho\nu}\left(S^{\rho}\otimes S^{\nu}\right).
$$

$$\chi_{\lambda}\downarrow^{S_n}_{S_k\times
S_{n-k}}=\sumlim_{\rho \vdash k, \nu \vdash n-k}{c_{\rho
\nu}^{\lambda}\chi_{(\rho,\nu)}}.$$

\ecl \qed

We use now the Frobenius reciprocity to obtain the multiplicity of
any irreducible representation of $S_n \times S_n$ in
$\alpha_{X_n^k}$.

\bpr Let $0\leq k\leq n$ and $\lambda,\mu\vdash n$.
\begin{multline}
m\big( (\lambda,\mu),\alpha_{X^k_n}\big)=\\
=\langle\chi_{\lambda}\downarrow^{S_n}_{S_k\times S_{n-k}},\chi_{\mu}\downarrow^{S_n}_{S_k\times S_{n-k}}\rangle=\\
=\sum_{\rho\vdash k,\nu\vdash n-k}c^{\lambda}_{\rho\nu}c^{\mu}_{\rho\nu}\nonumber
\end{multline}
\epr
\begin{proof}
\begin{multline}
m\left(\left(\lambda,\mu\right),\alpha_{X^k_n}\right)=\langle\chi_{\alpha_{X^k_n}},\chi_{(\lambda,\mu)}\rangle=\\
=\langle\chi_{\tilde{\omega}_{n,k}}\uparrow^{S_n\times S_n}_{(S_k\times S_{n-k})\times (S_k\times S_{n-k})},\chi_{(\lambda,\mu)}\rangle=\\
=\langle\chi_{\tilde{\omega}_{n,k}},\chi_{(\lambda,\mu)}\downarrow^{S_n\times S_n}_{(S_k\times S_{n-k})\times (S_k\times S_{n-k})}\rangle=\\
=\frac{1}{(k!)^2((n-k)!)^2}\sum_{(\pi_k\pi_{n-k},\sigma_k\sigma_{n-k})}\chi_{\tilde{\omega}_{n,k}}(\pi_k\pi_{n-k},\sigma_k\sigma_{n-k})\chi_{\lambda}(\pi_k\pi_{n-k})\chi_{\mu}(\sigma_k\sigma_{n-k})=\\
=\frac{1}{(k!)^2((n-k)!)^2}\sum_{\pi_k\pi_{n-k}\sim\sigma_k\sigma_{n-k}}\chi_{\tilde{\omega}_{n,k}}(\pi_k\pi_{n-k},\pi_k\pi_{n-k})\chi_{\lambda}(\pi_k\pi_{n-k})\chi_{\mu}(\pi_k\pi_{n-k})=\\
=\frac{1}{(k!)^2((n-k)!)^2}\sum_{\pi_k\pi_{n-k}}|C^{k\times (n-k)}(\pi_k\pi_{n-k})||C^k_{\pi_k}||C^{n-k}_{\pi_{n-k}}|\chi_{\lambda}(\pi_k\pi_{n-k})\chi_{\mu}(\pi_k\pi_{n-k})=\\
=\frac{1}{k!(n-k)!}\sum_{\pi_k\pi_{n-k}\in S_k\times S_{n-k}}\chi_{\lambda}(\pi_k\pi_{n-k})\chi_{\mu}(\pi_k\pi_{n-k})=\\
=\langle\chi_{\lambda}\downarrow^{S_n}_{S_k\times
S_{n-k}},\chi_{\mu}\downarrow^{S_n}_{S_k\times
S_{n-k}}\rangle=\sum_{\rho\vdash k,\nu\vdash
n-k}c^{\lambda}_{\rho\nu}c^{\mu}_{\rho\nu}\nonumber
\end{multline}
The last equality follows from Claim \ref{LR}.
\end{proof}

By the definition of $X_n^k$ we have
$\alpha_{B_n}=\bigoplus_{k=0}^n \alpha_{X^{k}_{n}}$ and thus:

\bco
$$
m\left(\left(\lambda,\mu\right),\alpha_{B_n}\right)=\sum_{k=0}^{n}\sum_{\rho\vdash k,\nu\vdash n-k}c^{\lambda}_{\rho\nu}c^{\mu}_{\rho\nu}.
$$
\eco \qed

There is a natural mapping between the sets $H_n^k$ and $X_n^k$
defined by:
$$
H_n^k\ni \pi U_{n,k}\sigma\stackrel{\tilde{T}}{\longmapsto}\pi \tilde{U}_{n,k}\sigma\in X_n^k
$$
One can verify that $\tilde{T}$ is well defined. Moreover,
$\tilde{T}$ commutes with the action $\alpha$ of $S_n\times S_n$
on $X^k_n$, i.e.:
$$
\tilde{T}(\pi A\sigma)=\pi \tilde{T}(A)\sigma\,\,\textrm{for any}\,\,A\in X^n_k.
$$

It is easy to see that $\tilde{T}$ is also surjective and thus it
induces epimorphisms of modules from the $S_n\times S_n$-module
$\alpha_{H_n^k}$ to the $S_n \times S_n$- module $\alpha_{X_n^k}$
and from the $S_n$-module $\beta_{H_n^k}$ to the $S_n$- module
$\beta_{X_n^k}$. Note also that for $k=0$ this mapping is the
identity mapping since $H_n^0=X_n^0=S_n$ and for $k=1$ this
mapping is bijective. We conclude: \bcl
$$
m\left(\left(\lambda,\mu\right),\alpha_{H_n^k}\right)\geq m\left(\left(\lambda,\mu\right),\alpha_{X_n^k}\right)
$$
\ecl \qed

This implies that if
$$
\sum_{\rho\vdash k,\nu\vdash n-k}c^{\lambda}_{\rho\nu}c^{\mu}_{\rho\nu}\neq 0
$$
then $|\lambda\bigtriangleup\mu|\leq 2k$. This can also be seen by
the combinatorial interpretation of the Littlewood-Richardson
coefficients.

\subsection{The action $\beta$ on colored permutations}
 Recall that  every matrix $B\in B_n$ can be written uniquely
in the form $B=Z\pi$ for some $\pi\in S_n$ and some $Z\in
 C_2^n$. There exists a natural epimorphism from $B_n$
onto $S_n$ defined by omitting the minuses:
$$
p\,:\,B_n\longrightarrow S_n \qquad\qquad p(Z\pi)=\pi.
$$
If we restrict $p$ to $X_n^k$ we obtain a surjective mapping from
$X_n^k$ onto $S_n$ which commutes with the action $\alpha$ of
$S_n\times S_n$ on $X^k_n$ (and clearly also commutes with the
action $\beta$ of $S_n$  on $X_n^k$ by conjugation). It gives us a
surjective homomorphism from the representation $\beta_{X_n^k}$
 onto the conjugacy representation
 representation of $S_n$ denoted by $\psi$. Therefore, using the result of \cite{F}, we
have
$$
m\left(\lambda,\beta_{X_n^k}\right)\geq m\left(\lambda,\psi\right)>0\quad\textrm{for any}\quad  \lambda\vdash n,
$$
$$
m\left(\lambda,\beta_{B_n}\right)=
\sum_{k=0}^{n}m\left(\lambda,\beta_{X_n^k}\right)>0\quad\textrm{for any}\quad  \lambda\vdash n.
$$
Although the calculation of $\chi_{\beta_{X_n^k}}$ is rather
involved the asymptotic results of [R] and [AF] can be generalized
for the representations $\beta_{X_n^k}$ and $\beta_{B_n}$. We
start by presenting the generalization of Theorem R1:

 \bpr\label{estimate_betaxnk}
 For any
$0<\varepsilon<1$ there exist $0<\delta(\varepsilon)$ and
$N(\varepsilon)$ such that, for any partition $\lambda$ of
$n>N(\varepsilon) $ with max$\{ {\lambda_1\over
n}\,,\,{\lambda_1'\over n}\}\leq\delta(\varepsilon)$,
$$
1-\varepsilon<\frac{m(\lambda,\beta_{X^k_n})}{{n\choose k}
f^\lambda}<1+\varepsilon\,\,,
$$
$$
1-\varepsilon<\frac{m(\lambda,\beta_{B_n})}{2^n
f^\lambda}<1+\varepsilon\,\,.
$$
\epr
\begin{proof}
Let $\pi \in S_n$. We have to estimate
$\chi_{\beta_{X_n^k}}(\pi)=\{B\in X_n^k\ \mid B \pi =\pi B\}$. If
$B \in X_n^k$ commutes with $\pi$ then $p(B)$ commutes with
$p(\pi)=\pi$ and thus $$|\{B \in X_n^k \mid \pi B=B \pi\}| \leq
|\{B \in X_n^k \mid p(B) \in C_{\pi}\}| = {n \choose
k}|C_{\pi}|\,\,.$$
By the same considerations we get:

$$
\chi_{\beta_{B_n}}(\pi)=\#\left\{B\in B_n\,|\,\pi
B=B\pi\right\}\leq 2^n|C_{\pi}|\,\,.
$$

In \cite{R} it is shown that for any $0<\varepsilon<1$ there exist
$0<\delta(\varepsilon)$ and $N(\varepsilon)$ such that, for any
partition $\lambda$ of $n>N(\varepsilon) $ with max$\{
{\lambda_1\over n}\,,\,{\lambda_1'\over
n}\}\leq\delta(\varepsilon)$,
$$
|m(\lambda)-f^{\lambda}|=|\sum_{C\in
\hat{S}_n}\chi_{\lambda}(C)-f^{\lambda}|=|\sum_{C\neq
id}\chi_{\lambda}(C)|\leq \varepsilon f^{\lambda}\,\,,
$$
which immediately implies Theorem R1. (Here $m(\lambda)$ denotes the multiplicity of the irreducible representation $\lambda$ in the conjugacy representation of $S_n$.) In our case we have $\chi_{\beta_{X_n^k}}(e)=|X^k_n|={n\choose k}n!$ and therefore
\begin{multline}
|m(\lambda,\beta_{X^k_n})- {n\choose k}f^{\lambda}|=|{1\over n!}\sum_{e\neq\pi\in S_n}\chi_{\lambda}(\pi)\chi_{\beta_{X_n^k}}(\pi)|\leq\\
\leq {n\choose k}{1\over n!}\sum_{e\neq\pi\in S_n}|\chi_{\lambda}(\pi)||C_{\pi}|={n\choose k}\sum_{\{e\}\neq C\in\hat{S}_n}|\chi_{\lambda}(C)|\leq {n\choose k}\varepsilon f^{\lambda}\nonumber
\end{multline}
which establishes the first statement. Substituting $2^n$ instead
of ${n\choose k}$ in this calculation we get the second statement.
\end{proof}

The generalization of Theorem \cite{AF} is as follows:

\bpr In the notations of Theorem \cite{AF}
$$
\lim_{n\to\infty}\frac{\|{n\choose k}\chi_R^{(n)}\|}{\|\chi_{\beta_{X^k_n}}\|}=\lim_{n\to\infty}\frac{\|2^n\chi_R^{(n)}\|}{\|\chi_{\beta_{B_n}}\|}=1\,,
$$
$$
\lim_{n\to\infty}\frac{\langle\chi_R^{(n)},\chi_{\beta_{X^k_n}}\rangle}{\|\chi_R^{(n)}\|\cdot\|\chi_{\beta_{X^k_n}}\|}=
\lim_{n\to\infty}\frac{\langle\chi_R^{(n)},\chi_{\beta_{B_n}}\rangle}{\|\chi_R^{(n)}\|\cdot\|\chi_{\beta_{B_n}}\|}=1\,,
$$
where $k$ is bounded or tends to infinity remaining less than $n$.
\epr
\begin{proof} Recall that
\begin{multline}
\|\chi_R^{(n)}\|^2=\langle\chi_R^{(n)},\chi_R^{(n)}\rangle=n!\,\,,\,\,\langle\chi_R^{(n)},\chi_{\beta_{X^k_n}}\rangle=\chi_{\beta_{X^k_n}}(e)=
{n\choose k}n!\\
\textrm{and}\,\,\langle\chi_R^{(n)},\chi_{\beta_{B_n}}\rangle=\chi_{\beta_{B_n}}(e)=2^n
n!\,.\nonumber
\end{multline}
Using the inequality $\chi_{\beta_{X_n^k}}(\pi)\leq {n\choose
k}|C_{\pi}|$ of Proposition \ref{estimate_betaxnk} we can write
$$
\|\chi_{\beta_{X_n^k}}\|^2={1\over n!}\sum_{\pi\in S_n}\left(\chi_{\beta_{X_n^k}}(\pi)\right)^2\leq {1\over n!}{n\choose k}^2\sum_{\pi\in S_n}|C_{\pi}|^2={n\choose k}^2 n!\sum_{C\in\hat{S}_n}{1\over |C|}.
$$
Also we obviously have
$$
\|\chi_{\beta_{X_n^k}}\|^2={1\over n!}\sum_{\pi\in S_n}\left(\chi_{\beta_{X_n^k}}(\pi)\right)^2\geq {1\over n!}\left(\chi_{\beta_{X_n^k}}(e)\right)^2={n\choose k}^2 n!\,.
$$
Taking the two last inequalities together we have:
$$
\frac{1}{\sqrt{\sum_{C\in\hat{S}_n}{1\over |C|}}}=\frac{{n\choose k}\sqrt{n!}}{{n\choose k}\sqrt{n!}\sqrt{\sum_{C\in\hat{S}_n}{1\over |C|}}}\leq \frac{\|{n\choose k}\chi_R^{(n)}\|}{\|\chi_{\beta_{X^k_n}}\|}\leq \frac{{n\choose k}\sqrt{n!}}{{n\choose k}\sqrt{n!}}=1\,.
$$
Also
$$
\frac{{n\choose k}n!}{\sqrt{n!}{n\choose k}\sqrt{n!}\sqrt{\sum_{C\in\hat{S}_n}{1\over |C|}}}\leq \frac{\langle\chi_R^{(n)},\chi_{\beta_{X^k_n}}\rangle}{\|\chi_R^{(n)}\|\cdot\|\chi_{\beta_{X^k_n}}\|}=\frac{{n\choose k}n!}{\sqrt{n!}\|\chi_{\beta_{X^k_n}}\|}\leq \frac{{n\choose k}n!}{\sqrt{n!}{n\choose k}\sqrt{n!}}=1\,.
$$
Substituting $2^n$ instead of ${n \choose k}$ we get similar
inequalities for $\|\chi_{\beta_{B_n}}\|$. In \cite{AF} it is
proved that
$$
\lim_{n\to\infty}\sum_{C\in\hat{S}_n}{1\over |C|}=1
$$
and thus using the above inequalities the proof is finished.
\end{proof}
These asymptotic results can be also obtained for the action
$\beta$ (conjugation by permutations) on the group
 $C_r \wr S_n$.
Similarly to $X_n^k\subset B_n$ define the sets $Y_n^k\subset C_r
\wr S_n$: \bde
$$
Y^{k}_{n}=\{A\in C_r\wr S_n\,|\,A \,\,\text{has  exactly
}\,k\,\,\,\text{entries}\neq 0\,,\,1\}.
$$
\ede

Note that the sets $Y^k_n$ form a partition of $C_r \wr S_n$ and
$Y_n^k =n!{n\choose k}(r-1)^k$.  The sets $Y^k_n$ are closed under
the action $\alpha$ of $S_n\times S_n$ but they are not transitive
under this action.

Consider $C_r^n$ as the group of diagonal matrices with the
entries of the form $\omega^{\ell}$ (where $\omega=\exp{{2\pi
i\over r}}$ -- the primitive $r$-th root of unity and $0\leq \ell
<r$) on the diagonal. Then each matrix $A\in C_r\wr S_n$ can be
uniquely written as $A=Z\sigma $ for some $\sigma\in S_n$ and some
$Z\in C_r^n$. Just as in the case of $B_n$, we consider the
epimorphism $p\,:\,B_n\longrightarrow S_n$ defined by:
$p(Z\sigma)=\sigma$.

$p$ induces an epimorphism of modules between $\beta_{Y_n^k}$ and
the conjugacy representation of $S_n$.

We conclude, using the result of \cite{F}:
$$
m\left(\lambda,\beta_{Y_n^k}\right)\geq m\left(\lambda,\psi\right)>0\quad\textrm{for any}\quad  \lambda\vdash n,
$$
$$
m\left(\lambda,\beta_{C_r \wr
S_n}\right)=\sum_{k=0}^{n}m\left(\lambda,\beta_{Y_n^k}\right)>0\quad\textrm{for
any}\quad  \lambda\vdash n.
$$

The Theorems R1 and AF are obtained in a way similar to the one we
used for $B_n$: \bpr In the conditions and notations of Theorem R1
$$
1-\varepsilon<\frac{m(\lambda,\beta_{Y^k_n})}{{n\choose k}(r-1)^k
f^\lambda}<1+\varepsilon,
$$
$$
1-\varepsilon<\frac{m(\lambda,\beta_{C_r \wr S_n})}{r^n
f^\lambda}<1+\varepsilon.
$$
In the notations of Theorem AF
$$
\lim_{n\to\infty}\frac{\|{n\choose
k}(r-1)^k\chi_R^{(n)}\|}{\|\chi_{\beta_{Y^k_n}}\|}=\lim_{n\to\infty}\frac{\|r^n\chi_R^{(n)}\|}{\|\chi_{\beta_{C_r
\wr S_n}}\|}=1\,,
$$
$$
\lim_{n\to\infty}\frac{\langle\chi_R^{(n)},\chi_{\beta_{Y^k_n}}\rangle}{\|\chi_R^{(n)}\|\cdot\|\chi_{\beta_{Y^k_n}}\|}=
\lim_{n\to\infty}\frac{\langle\chi_R^{(n)},\chi_{\beta_{C_r\wr
S_n}}\rangle}{\|\chi_R^{(n)}\|\cdot\|\chi_{\beta_{C_r \wr
S_n}}\|}=1\,,
$$
where $k$ is bounded or tends to infinity remaining less than $n$.
\epr
\section{Appendix}

We defined the sets $H_n^k$ in order to study the representations $\beta_{H_n^k}$ and $\alpha_{H_n^k}$.
These representations can be obtained from the action $\beta$ of $S_n$ and the action $\alpha$ of $S_n\times S_n$ on some
other sets of matrices instead of $H_n^k$:\\
\\
1) In the definition of the matrix $U_{n,k}$ change only one thing: put the upper
 left $k\times (n-k)$ block to be a zero matrix instead of ones in $U_{n,k}$.
 Taking the orbit of this matrix under the action (1) of $S_n\times S_n$ we have
  a natural bijection between the obtained set (this orbit) and  $H_n^k$ which preserves
   the action (1). So we get the representations which are isomorphic to the representations
   $\beta_{H_n^k}$ and $\alpha_{H_n^k}$.\\
\\
2) Instead of $H_n^k$ we can take the following subset of $C_{k+1}\wr S_n$:\\
\begin{multline}
\left\{A\in C_{k+1}\wr S_n\,|\, A\, \textrm{ has $\omega,\omega^2,\ldots,\omega^k$ exactly once}\right.\\
 \left. \textrm{and all the other is as in a permutation matix.}\Big\}\right.\nonumber
\end{multline}
 The representations $\alpha$ and $\beta$ obtained from the corresponding actions on these sets are also isomorphic to $\alpha_{H_n^k}$ and $\beta_{H_n^k}$.


\noindent
\begin{thebibliography}{FLM}
\newcommand{\au}{\sc}
\newcommand{\ti}{\it}
\bibitem[AF]{AF}
{\au R.\ M.\ Adin and A.\ Frumkin},
{\ti The conjugacy character tends to be regular},
Israel J.\ Math.~{\bf59} (1987), 234--240.
\bibitem[CS]{CS}
{\au Y.\ Cherniavsky and M.\ Sklarz},
{\ti On conjugation action of $S_n$ on invertible matrices}, Preprint, available on the website:
www.math.biu.ac.il/\~{}cherniy
\bibitem[F]{F}
{\au A.\ Frumkin},
{\ti Theorem about the conjugacy representation of $S_n$},
Israel J.\ Math.~{\bf55} (1986), 121--128.
\bibitem[I]{I}{\au Iwahori N.}{\ti On an equivalence relatioion
among (0,1)-matrices}, Sci.Papers College Gen.Ed.Univ.Tokyo ~{\bf
21}(1971), 1-9.
\bibitem[Li]{Li}{\au Li J.S.}{\ti Permutation equivalence of
(0,1)- square matrices}, Acta Math.Sinica  ~{\bf 26} (1983), no.5,
586-596.
\bibitem[Md]{Md}
{\au I.\ G.\ Macdonald},
{\ti Symmetric Functions and Hall Polynomials},
second edition, Oxford Math.\ Monographs,
Oxford Univ.\ Press, Oxford, 1995.

\bibitem[MM]{MM}{\au Mader A., Mutzbauer O.},{\ti Double orderings
of (0,1)-matrices}, Ars Combin. ~{\bf 61}(2001),81-95.

\bibitem [R]{R}
{\au Y.\ Roichman},
{\ti Decomposition of the conjugacy representation of the symmetric groups},
Israel J.\ Math.~{\bf97} (1997), 305--316.
\bibitem[Sa]{Sag}
 B.\ Sagan, {\it The symmetric group: representations,
combinatorial algorithms, and symmetric functions}, Springer
Verlag (2001).

\end{thebibliography}
\end{document}